\documentclass[11pt]{article}

\usepackage{amsfonts}

\newcommand{\CR}{\hbox{{$\cal R$}}} 

\newcommand{\R}{\mathbb{R}}

\newcommand{\C}{\mathbb{C}}

\newcommand{\Z}{\mathbb{Z}}

\renewcommand{\H}{\mathbb{H}}

\newcommand{\del}{\partial}

\newcommand{\dsl}{{\del\!\!\!/}}

\newcommand{\h}{{\scriptstyle\frac{1}{2}}}

\newcommand{\isom}{{\cong}}

\newcommand{\eps}{{\epsilon}}

\newcommand{\tens}{\mathop{\otimes}}

\newcommand{\ad}{{\rm ad}}

\newcommand{\lcross}{{>\!\!\!\triangleleft}}

\newcommand{\proof}{{\bf Proof\ }}

\newcommand{\eproof}{$\quad \diamond$\bigskip}

\newcommand{\und}[1]{{\underline {#1}}}

\newcommand{\eqn}[2]{\begin{equation}#2\label{#1}\end{equation}}

\newtheorem{lemma}{Lemma}[section]

\newtheorem{propos}[lemma]{Proposition}

\newtheorem{corol}[lemma]{Corollary}

\begin{document}

\hsize 17truecm \vsize 24truecm \font\twelve=cmbx10 at 13pt

\font\eightrm=cmr8 \baselineskip 18pt

{\ }\qquad \hskip 4.3in \vspace{.2in}
\begin{center} {\LARGE Clifford Algebras Obtained by Twisting
of Group Algebras}\\ \baselineskip 13pt{\ }
{\ }\\ Helena Albuquerque\footnote{Supported by CMUC-FCT} \\{\ }\\
Departamento de
Matematica-Faculdade de Ciencias e Tecnologia\\ Universidade de
Coimbra, Apartado 3008\\ 3000 Coimbra, Portugal\\ lena@mat.uc.pt{\ }\\
+\\{\ }\\ Shahn Majid \\ {\
}\\ School of Mathematical Sciences\\ Queen Mary, University of
London\\ Mile End Rd, London E1 4NS, UK\\
www.maths.qmw.ac.uk/$\tilde{\ }$majid\end{center}
\begin{center}
October 2000
\end{center}

\begin{quote}\baselineskip 13pt
\noindent{\bf Abstract} We investigate the construction and properties
of Clifford algebras by a similar manner as our previous
construction of the octonions, namely as a twisting of group
algebras of $\Z_2^n$ by a cocycle. Our approach is more general
than the usual one based on generators and relations. We obtain
in particular the periodicity properties and a new construction of
spinors in terms of left and right multiplication in the
Clifford algebra. \end{quote}

\section{Introduction}

In \cite{AlbMa:qua} we have constructed the octonions and other
Cayley algebras by a twisting procedure applied to group algebras
$k G$ by a 2-cochain $F$ on the group. The failure of the cochain to
be a group cocycle controls the nonassociativity of these algebras.
On the other hand, in the case when the
cochain is actually a group cocycle the associativity will be
preserved by twisting. It has already been observed that, in particular, if
one uses cocycles which are the quadratic part of the cochains in
\cite{AlbMa:qua} that define the Cayley algebras then one has in fact their
associated Clifford algebras.

In this paper we will explore this construction further, using
this point of view to give a new derivation of known results about Clifford
algebras and to generalise them. Results include the periodicity
theorems for Clifford algebras and a novel construction for their
spinor representations. The paper begins in Section~2 by showing how
a specific cocycle gives the usual Clifford algebras and many of their
properties in a more direct manner. Section~3
contains some constructions for general groups and cocycles that
include and generalise the periodicity theorems. Section~4 contains
spinor constructions for usual Clifford algebras obtained by our
methods. As an example, spinors in 4 dimensions are naturally
described as quaternion valued functions.

\subsection{Preliminaries}

Let $k$ be a field with characteristic not 2. Let ${\bf q}$ be a
nondegenerate quadratic form on a vector space $V$ over $k$ of
dimension $n$. It is known\cite{Lam:quad} that there is an orthogonal basis
$\{e_1,\cdots,e_n\}$, say, of $V$ such that ${\bf q}(e_i)=q_i$
for some $q_i\ne 0$. The Clifford algebra $C(V,{\bf q})$,
\cite{Wen:cliff},
is the associative algebra generated by $1$ and $\{e_i\}$ with the
relations
\[ e_i^2=q_i.1,\quad e_ie_j+e_je_i=0,\quad \forall i\ne j\]
We identify $k$ and $V$ inside $C(V,{\bf q})$ in the obvious way.
The dimension of $C(V,{\bf q})$ is $2^n$ and it has a canonical
basis
\[ \{e_{i_1}\cdots e_{i_p}|\, 1 \le\, i_1 <  i_2\cdots <i_p\le n\}.\]

If we assume that $q_i=\pm 1$ (without loss of generality over $\R$,
for example) then for $n=1$ we have two cases: If $q_1=-1$ we have
$C(k,q_1)$ the
algebraic complex numbers where we adjoin $i=e_1$ with relation
$i^2=-1$. If $q_1=1$ we have the group
algebra of $\Z_2$ with $e_1$ the   nontrivial element. Setting
$e_\pm=(1\pm e_1)/2$ we have equivalently two projections
\[ e_\pm^2=e_\pm,\quad e_+e_-=e_-e_+=0\]
and $C(k,q_1)= k\oplus k$ (the hyperbolic complex numbers
over $\R$).

If $n=2$ we have
\[ e_1e_2=-e_2e_1,\quad e_i^2=q_i,\quad (e_1e_2)^2=-q_1q_2,\quad
(e_1e_2)e_1=-q_1e_2,\quad (e_1e_2)e_2=q_2e_1.\] Hence for
$q_1=q_2=-1$ we have the $k$-algebra of quaternions $\H$ with
$i=e_1$, $j=e_2$ and $k=e_1e_2$. If $q_1=q_2=1$ we have the
matrix algebra $M_2(k)$ with $e_1=\left(\begin{matrix}{1 & 0\cr
0&-1}\end{matrix}\right)$ and $e_2=\left(\begin{matrix}{0&1\cr
1&0}\end{matrix}\right)$. If $q_1=1$ and $q_2=-1$ we have $M_2(k)$
similarly with $e_2=\left(\begin{matrix}{0&1\cr
-1&0}\end{matrix}\right)$ instead. Thus
\[ C(0,2)\isom\H,\quad C(2,0)\isom C(1,1)\isom M_2(k)\]
where $C(r,s)$ denotes the algebra with $r$ of the $\{q_i\}$ equal
to $+1$ and $s$ equal to $-1$.

\section{$C(V)$ as twisting}

{}From \cite{AlbMa:qua} we recall that if $G$ is a group and
$F:G\times G\to k$ a nowhere-zero function with $F(e,x)=F(x,e)=1$
(a cochain), where $e$ is the group identity, then we define a new
algebra $k_FG$ which has the same vector space as the group
algebra $kG$ (namely with basis labelled by $G$) but a different
product, namely \eqn{prod}{ a\cdot b= F(a,b)ab,\quad \forall
a,b\in G.} This arises naturally as an algebra in the category of
comodules of the dual-quasiHopf algebra $(kG,\del F)$. When $F$
is a cocycle this gives a cotriangular usual Hopf algebra
structure on $kG$. We first construct this.

\begin{propos} Let $G=\Z_2^n$. There is a 2-cocycle $F\in Z^2(G,k)$
defined by
\[ F(x,y)=(-1)^{\sum_{j<i}x_iy_j}
\prod_{i=1}^nq_i^{x_iy_i}\] where $x=(x_1,\cdots x_n)\in \Z_2^n$
and twists $k G$ into a cotriangular Hopf algebra\cite{Ma:book} with
cotriangular structure
\[ \CR(x,y)=(-1)^{\rho(x)\rho(y)+x\cdot
y},\] where $\rho(x)=\sum_i x_i\in\Z$ and $x\cdot y$ is the dot
product of $\Z_2$-valued vectors.
\end{propos}
\proof  We then define the cochain $F$ as shown which is
manifestly invertible as $q_i\ne 0$ and the identity when either
argument is zero. Using the notation $x=(x_1,\cdots,x_n)$, it is
clear that
\[ \del F(x,y,z)=\frac{F(x,y)F(x+y,z)}{F(y,z)F(x,y+z)}=\prod_{i=1}^n
q_i^{x_iy_i+(x_i+y_i)z_i -y_iz_i-x_i(y_i+z_i)}=1\] since the
$(-1)$ factors certainly do not contribute by linearity. The
remaining factors again cancel because the exponents are
controlled by a bilinear form. Hence by the twisting theory for
Hopf algebras\cite{Dri:qua}\cite{Ma:book} the group algebra $k G$
becomes a
cotriangular Hopf algebra with the same Hopf algebra structure
(namely every element of the group has diagonal coproduct) but
with $\CR: k G\tens k G\to k$ defined by $\CR(x,y)={F(x,y)\over
F(y,x)}$ on the basis elements of $k G$ (labeled by $G$). This
readily comes out as shown. \eproof

The modification of product (\ref{prod}) from $kG$ to a new
algebra $k_F G$ is an application of twisting to comodule
algebras\cite{GurMa:bra} and ensures that all its structure maps are
morphisms in the symmetric monoidal category\cite{Mac:cat} of
$(kG,\del F)$-comodules.
This category is the same as that of $G$-graded spaces equipped
with a generalised transposition given by $\CR$ and associativity
by $\del F$. Because in our case the coboundary $\del F=1$, the
new algebra remains associative in the usual sense.

Applying this construction, the original product of basis
elements in $k\Z_2^n$ corresponds to the addition in $\Z_2^n$. We
now define $k_F \Z_2^n$ as the same vector space as $k \Z_2^n$
but with the new product modified by $F$. For clarity we denote
the basis elements of $k \Z_2^n$ by
\eqn{ex}{ e_x=e_1^{x_1}\cdots e_n^{x_n}}
where the unmodified algebra $k \Z_2^n$ is generated by $e_i$
mutually commuting and with $e_i^2=1$.

\begin{propos}
The algebra $k_F\Z_2^n$ can be identified with $C(V,{\bf q})$,
i.e. the latter is an algebra in the symmetric monoidal category
of $\Z_2^n$-graded spaces defined by $\CR$.
\end{propos}
\proof We identify the elements $e_i$ in the two algebras and
note that $e_i\cdot e_j=e_ie_j$ if $i<j$ since $F(x,y)=1$ in this
case (where $x_i=1=y_j$ and other entries are zero). Hence we can
identify also the basis elements
\[ e_1^{x_1}\cdot e_2^{x_2}\cdots e_n^{x_n}=e_x\]
of the two algebras. After that it remains only to check that the
products coincide for other products, which can be check
inductively from the generators. Here one has
\[ e_i\cdot e_j=-e_ie_j=-e_je_i=-e_j\cdot e_i,\quad e_i\cdot
e_i=q_ie_i^2=q_i 1,\quad \forall i, j<i.\] \eproof

We can now apply some of the results in \cite{AlbMa:qua} albeit in
the associative setting since $\phi=\del F=1$.

\begin{corol}
The algebra $C(V,{\bf q})$ is commutative under the generalised
transposition or `braiding' of the symmetric monoidal category
defined by $\CR$, i.e. $e_x\cdot e_y=\CR(x,y)e_y\cdot e_x$. In
particular,
\[ e_x \cdot e_y=\cases{e_y\cdot e_x &{\rm if}\  $x\cdot y=\rho(x)\rho(y)$
\ {\rm (mod 2)}\cr -e_y\cdot e_x & {\rm else.}\cr}\]
\end{corol}
\proof This is an immediate consequence of the above construction
in terms of $F$ since $a\cdot b=F(a,b)ab=F(a,b)ba={F(a,b)\over
F(b,a)}b\cdot a$. In our case it takes the form shown. The result
is also clearly true on repeated use of the anticommutation
relations of $C(V,{\bf q})$ but in our approach these  are
encoded concisely in $\CR$ as stated. Note that either $x\cdot
y=\rho(x)\rho(y)$ or $x\cdot y=\rho(x)\rho(y)+1$ mod 2, so two
basis elements either commute or anticommute. \eproof

We therefore apply various results about $k_FG$ algebras in
\cite{AlbMa:qua} to Clifford algebras. For example,

\begin{propos} $\CR$ in Proposition~2.1 is a group coboundary,
$\CR=\del \theta$ for cochain
\[ \theta(x)=(-1)^{\h \rho(x)(\rho(x)-1)}\]
with $\theta^2=1$. Hence by \cite[Lem. 3.4]{AlbMa:qua}
$\Theta(e_x)=\theta(x)e_x$
is a diagonal anti-involution on the Clifford algebra $C(V,{\bf
q})$. Explicitly, $\Theta$ is the order-reversal operation
\[ \Theta(e_1^{x_1}\cdots e_n^{x_n})=e_n^{x_n}\cdots e_1^{x_1}.\]
\end{propos}
\proof We show that $\CR(x,y)=\theta(x)\theta(y)/\theta(x+y)$,
where we write the group structure of $\Z_2^n$ additively. Viewing
$\rho(x)=\sum_i x_i$ in $\Z$, \eqn{rhoplus}{ \rho(x+y)=\sum_i
x_i(1-y_i)+(1-x_i)y_i=\rho(x)+\rho(y)-2 x\cdot y} and hence
\[ \theta(x+y)=(-1)^{\h\left(\rho(x+y)^2-\rho(x+y)\right)}=
(-1)^{\h\left(\rho(x)^2+\rho(y)^2+2\rho(x)\rho(y)-\rho(x)
-\rho(y)+2x\cdot
y\right)}=\theta(x)\theta(y)(-1)^{\rho(x)\rho(y)+x\cdot y}.\]
Hence $\theta$ defines an anti-involution $\Theta$. In $C(V,{\bf q})$
we can identify it with order-reversal on noting that
\[ e_n^{x_n}\cdots e_1^{x_1}=e_x (-1)^{\sum_{j<i}x_ix_j}=e_x
(-1)^{\h\rho(x)(\rho(x)-1)},\]where $\rho(x)^2=\sum_{i,j}
x_ix_j=2\sum_{j<i}x_ix_j+\sum_i x_i^2= 2\sum_{j<i}x_ix_j+\rho(x)$
since $x_i=0,1$. \eproof

It is clear that $\Theta$ is an isomorphism between $C(V,{\bf q})$ and
its opposite algebra, which in turn is of form $k_{F^{\rm
op}}\Z_2^n$ for $F^{\rm op}(x,y)=F(y,x)$. We also obtain from
$\rho$ the obvious $\Z_2$-grading  of Clifford algebras provided
by an order 2 automorphism. Here degree zero is the eigenspace
with eigenvalue 1 under the involution.

\begin{corol} $\sigma(e_x)=(-1)^{\rho(x)}e_x$ extended linearly is an
automorphism  of $C(V,{\bf q})$ in the form above and makes it
into a super-algebra with $\rho$ the super degree. If $n$ is even
then $\sigma$ is inner, being implemented by
$e_{(1,\cdots,1)}=e_1\cdots e_{n}$.  \end{corol} \proof The first
part is again an immediate consequence of the above construction
in terms of $F$ since $k_F\Z_2^n$ is $\Z_2^n$-covariant (it is an
algebra in the category of $\Z_2^n$-graded spaces). The map
$\rho:\Z_2^n\to \Z$ given by $\rho(x)=\sum_{i=1}^nx_i$ is
additive mod 2 (a group homomorphism to $\Z_2$) and therefore induces a
functor from the category of $\Z_2^n$-graded spaces to that of
$\Z_2$-graded ones. Under this any $\Z_2^n$-graded algebra is also
a $\Z_2$-graded or super algebra. The second part is a well-known
from generators and relations. In our case it comes about as
\[ e_{(1,\cdots,1)}^{-1}e_x e_{(1,\cdots,1)}
=e_{(1,\cdots,1)}^{-1}e_{(1,\cdots,1)}e_x (-1)^{\rho(x)
\rho(1,\cdots,1)+x\cdot(1,\cdots,1)}=e_x (-1)^{(n+1)\rho(x)}\]
using the braided-commutativity in Corollary~2.3. Note also that
\eqn{esq}{ e_{(1,\cdots,1)}^2=F((1,\cdots 1),(1,\cdots,1))
=(-1)^{n(n-1)\over 2}\prod_iq_i.}
\eproof

In terms of this superalgebra structure one can say that the
natural braiding on $C(V,{\bf q})$ defined by $\CR$ in Proposition~2.1
(with respect to which the Clifford algebra is braided-commutative) is
of the form \eqn{Psi}{ \Psi(e_x\tens e_y)=\Psi_{\rm
super}(e_x\tens e_y)(-1)^{x\cdot y}} where $\Psi_{\rm super}$
refers to the usual bose-fermi statistics or supertransposition.
There are of course many other applications of the superalgebra
structure.

\begin{corol} $C(V\oplus W,{\bf q}\oplus {\bf p})\isom
C(V,{\bf q})\und\tens C(W,{\bf p})$ as super algebras.
\end{corol}
\proof This is well-known from the point of view of generators
and relations. In our description it is clear from the form of
$F$ in Proposition~2.1 as
\[ F((x,x'),(y,y'))=F(x,y)F(x',y')(-1)^{\rho(x')\rho(y)}\]
where $\{e_{x}\}$ is a basis of $V$ and $\{e_{x'}\}$ of $W$,
say.  Hence the algebra product has the form $(a\tens c)(b\tens
d)=a\cdot b\tens c\cdot d (-1)^{\rho(c)\rho(b)}$ for the super
tensor product of super algebras. The notation ${\bf q}\oplus{\bf
p}$ indicates zero inner product between $V$ and $W$. \eproof

Finally, our approach also gives more explicit formulae for the
adjoint action and Pin groups. First of all our explicit form of
the product means that all the basis elements $e_x$ are
invertible in $C(V,{\bf q})$, as are generic linear combinations
for $q_i=\pm 1$. The latter fact is because the products all have
coefficients $\pm 1$ coming from the values of $F$. Clearly
\eqn{inv}{e_x^{-1}={e_x\over F(x,x)}=e_x
{(-1)^{\sum_{j<i}x_ix_j}\over \prod_{i=1}^nq_i^{x_i}}= {\Theta(e_x
)\over \prod_{i=1}^nq_i^{x_i}}} in terms of the anti-involution
above. We recall that the Clifford group of $V$ consists of the
invertible elements of $C(V,{\bf q})$ that leave $V$ stable under
the adjoint action.

\begin{corol} The adjoint action defined by
$\ad_{a}(e_y)=\sigma(a)e_y a^{-1}$ for all invertible $a\in
C(V,{\bf q})$ takes the explicit form \[
\ad_{e_x}(e_y)=(-1)^{\rho(x)(\rho(y)+1)}(-1)^{x\cdot y} e_y.\]
\end{corol} \proof This follows immediately from the
braided-commutativity, i.e. from the form of $\CR$ in Proposition~2.1
and Corollary~2.3. Thus $(-1)^{\rho(x)}e_x\cdot e_y\cdot
e_x^{-1}=(-1)^{\rho(x)}(-1)^{\rho(x)\rho(y)+x\cdot y}e_y\cdot
e_x\cdot e_x^{-1}$. \eproof

Also, using $S$, one defines $\lambda:C(V,{\bf q})\to k$ by
$e_x\sigma\circ \Theta e_x=\lambda(e_x)1$. From our formula for
inverses we obtain explicitly \eqn{lam}{
\lambda(e_x)=(-1)^{\rho(x)}\prod_{i=1}^nq_i^{x_i}.} By definition
the group ${\rm Pin}(V)$ is the subgroup of the Clifford group
with $\lambda=\pm 1$ and clearly includes all the $e_x$ when
$q_i=\pm 1$. The even part of this is the spin group. These
groups map surjectively onto $O(V)$ and $SO(V)$ via $\ad$.

\section{Clifford process}

We now use the above convenient description of Clifford algebras
to express a `doubling process' similar to the Dickson process for
division algebras. Thus, let $A$ be a finite-dimensional algebra
with identity $1$ and $\sigma$ an involutive automorphism of $A$.
For any fixed element $q\in k^*$ there is a new algebra of twice
the dimension,
\[ \bar A=A\oplus Av,\quad (a+bv)\cdot(c+dv)=a\cdot
c+q b\cdot \sigma(d)+(a\cdot d+b\cdot \sigma(c))v\] with a new
involutive automorphism
\[ \bar\sigma(a+vb)=\sigma(a)-\sigma(b)v.\]
We will say that $\bar A$ is obtained from $A$ by {\em Clifford
process}, see \cite{Wen:cliff}. We consider this initially for not
necessarily associative algebras and then find conditions for
associativity to be preserved.

\begin{propos} Let $G$ be a finite Abelian group and $F$ a
cochain as above, so $k_FG$ is a $G$-graded quasialgebra. For any
$s:G\to k^*$ with $s(e)=1$ and any $q\in k^*$, define $\bar
G=G\times\Z_2$ and
\[ \bar F(x,yv)=F(x,y)=\bar F(x,y),\quad \bar
F(xv,y)=s(y)F(x,y)\]
\[ \bar F(xv,yv)=q s(y) F(x,y),\quad \bar s(x)=s(x),\quad
\bar s(xv)=-s(x)\] for all $x,y\in G$. Here $x\equiv (x,e)$ and
$xv\equiv(x,\eta)$ where $\eta$ with $\eta^2=e$ is the generator
of the $\Z_2$. If $\sigma(e_x)=s(x)e_x$ is an involutive
automorphism then $k_{\bar F}\bar G$ is the Clifford process
applied to $k_FG$.
\end{propos}
\proof We clearly have a new cochain since $\bar
F(e,xv)=F(e,x)=1$ and $\bar F(xv,e)=s(e)F(x,e)=1$. The formulae
are fixed by reproducing the product of $\bar A$ in the
involutive case. Thus $F(x,yv)xyv=x\cdot yv=(x\cdot
y)v=F(x,y)xyv$, $F(xv,y)xvy=xv\cdot y=(x\cdot \sigma(y))v=s(y)(x\cdot
y)v=s(y)F(x,y)xyv$, etc.\eproof

It is easy to see that the special case where $s$ defines an involutive
automorphism on $k_FG$ is precisely the one where $s:G\to k^*$ is a
character with $s^2=1$.

\begin{propos}  For any $s:G\to k^*$ and $q\in k^*$ as above the $k_{\bar
F}\bar G$ given by the generalised Clifford process has
associator and braiding \[
\bar\phi(x,yv,z)=\bar\phi(x,y,zv)=\bar\phi(x,yv,zv)=\phi(x,y
,z)\]
\[ \bar\phi(xv,y,z)=\bar\phi(xv,yv,z)=\phi(xv,y,zv)= \phi(xv,yv,zv)=
\phi(x,y,z){s(yz)\over s(y)s(z)}\]
\[ \bar{\CR}(x,y)=\CR(x,y),\quad \bar{\CR}(xv,y)=s(y)\CR(x,y),\quad
\bar{\CR}(x,yv)={\CR(x,y)\over s(x)},\quad
\bar{\CR}(xv,yv)=\CR(x,y){s(y)\over s(x)}.\]
\end{propos}
\proof This is an elementary computation from the definitions of
$\phi,\CR$ for $k_FG$ and $k_{\bar F}\bar G$ and the form of
$\bar F$ above. For example,
\[ \bar\phi(xv,y,z)={\bar F(xv,y)\bar F(xvy,z)\over\bar F(y,z)\bar F(xv,y
z)}={s(y)F(x,y)s(z)F(xy,z)\over
F(y,z)s(yz)F(x,yz)}=\phi(x,y,z){s(y)s(z)\over s(yz)}\] as
stated. Similarly $\bar {\CR}(xv,yv)={\bar F(xv,yv)\over\bar
F(yv,xv)}={s(y)F(x,y)\over s(x)F(y,x)}={\CR(x,y)s(y)\over s(x)}$,
etc. \eproof

The merit of our approach is that these computations of the
associator and braiding are elementary but the properties of
$k_{\bar F}\bar G$ can be read off in terms of them. Thus it is
immediate that,

\begin{corol} If $s$ defines an involutive automorphism $\sigma$ then
$\bar\phi=1$ iff $\phi=1$, i.e. $k_{\bar F}\bar G$ is associative
iff $k_FG$ is.
\end{corol}
\proof In this case $\bar\phi$ and $\phi$ are given by the same
expressions independently of the placement of $v$. \eproof

Similarly, we gave in [1] conditions for $k_FG$ to be
alternative in terms of $\CR,\phi$. Using these, we have

\begin{corol} If $s$ defines an involutive automorphism $\sigma$ then
$k_{\bar F}\bar
G$
is alternative iff

(i) $k_FG$ is alternative

(ii) For all $x,y,z\in G$, either $\phi(x,y,z)=1$ or
$s(x)=s(y)=s(z)=1$.
\end{corol} \proof Since $\bar\phi,\bar{\CR}$ restrict to
$\phi,\CR$ it is immediate that $k_{\bar F}\bar G$ alternative
implies $k_FG$ alternative. Here alternativity of $k_FG$ is
explicitly the condition \[
\phi(x,y,z)+\CR(z,y)\phi(x,z,y)=1+\CR(z,y)\] \[
\phi^{-1}(x,y,z)+\CR(y,x)\phi^{-1}(y,x,z)=1+\CR(y,x)\] while for
$k_{\bar F}\bar G$ we have these and other cases such as
\[ \bar\phi(x,y,zv)+\bar{\CR}(zv,y)\bar\phi(x,zv,y)
=1+\bar{\CR}(zv,y)\]
or, from the above results,
\[ \phi(x,y,z)+s(y)\CR(z,y)\phi(x,z,y)=1+s(y)\CR(z,y).\]
Comparing, we see that in this case $(\phi(x,z,y)-1)(s(y)-1)=0$.
Similarly the content of the other cases of the condition for
$k_{\bar F}\bar G$ alternative is precisely that
$(\phi(x,z,y)-1)(s(x)-1)=0$ and $(\phi(x,z,y)-1)(s(z)-1)=0$ as
well.  Thus $k_{\bar F}\bar G$ is alternative iff $k_FG$ is {\em
and} for all $x,y,z$,
\[\phi(x,y,z)=1,\quad{\rm or}\quad s(x)=s(y)=s(z)=1.\]
\eproof

What this means is that either $k_FG$ is associative for the
conclusion to hold or, if not, then $\sigma$ has to be nontrivial
for some of the elements whose product fails to associate.

We now look at the case where $F,s,q$ are of the form
\eqn{Fz2}{ F(x,y)=(-1)^{f(x,y)},\quad s(x)=(-1)^{\xi(x)},\quad
q=(-1)^\eps}
for some $\Z_2$-valued functions $f,\xi$ and $\eps\in \Z_2$. We
also suppose that $G=\Z_2^n$ and use a vector notation.

\begin{lemma} For $G,F,s$ based on $\Z_2$, the generalised Clifford
process  yields the same form with $G=\Z_2^{n+1}$ and \[ \bar
f((x,x_{n+1}),(y,y_{n+1}))=f(x,y)+(y_{n+1}\eps+\xi(y))x_{n+1},\quad
\bar\xi(x,x_{n+1})=\xi(x)+x_{n+1}.\] \end{lemma} \proof Clearly
$\bar f(xv,yv)=\eps+\xi(y)+f(x,y)$ is the case where
$x_{n+1}=y_{n+1}=1$, while $\bar f(xv,y)= \xi(y)+f(x,y)$  is the
case where $x_{n+1}=1$ and $y_{n+1}=0$. The other two cases
require to yield $f(x,y)$. The four cases can then be expressed
together as stated using the field $\Z_2$. Similarly for
$\bar\xi$. \eproof

\begin{corol} Starting with $k$ and iterating the Clifford process
with
a choice of $q_i=(-1)^{\eps_i}$ at each step, we arrive at the
standard $C(V,{\bf q})$ in Proposition~2.1 and the standard
automorphism $\sigma(e_x)=(-1)^{\rho(x)}e_x$. \end{corol} \proof
We start with $f=0$ and $\xi=0$. Clearly $\xi(x)=\rho(x)$ after
$n$ steps independently of $\eps_i$. We also build up the
expression $\sum_{j<i}x_iy_j$ in $f$ as required, and additional
contribution to $f$ which gives the product the expression in
Proposition~2.1.
\eproof

Equivalently we can read Lemma~3.5 inductively. If we use the
notation $C(r,s)$ for the number of $\pm$ in the quadratic form
then,

\begin{corol} Starting with $C(r,s)$ the Clifford process with $q=1$ yields
$C(r+1,s)$. With $q=-1$ it gives $C(r,s+1)$. Hence any $C(m,n)$
with $m\ge r, n\ge s$ can be obtained from successive
applications of the Clifford process from $C(r,s)$. \end{corol}
\proof $\bar f$ in the lemma above, given that
$\xi(x)=\rho(x)$, is manifestly of the required form. Here
$\xi(y)x_{n+1}=\sum_{j<n+1}x_{n+1}y_j$ and $\eps y_{n+1}x_{n+1}$ gives
the extra factor in the product in the expression in Proposition~2.1.
\eproof

Note also that the definition of $\bar A$ can be written equally well
as some kind of `tensor product' $\bar A=A\tens_\sigma C(k,q)$
where $C(k, q)=k[v]$ with the relation $v^2=q$, and
$\tens_\sigma$ denotes that $\bar A$ factorises into these
subalgebras with the cross relations $va=\sigma(a)v$ for all
$a\in A$.  As $k_FG$
algebras we do not need to assume an involutive automorphism
$\sigma$ and clearly have $k_{\bar F}\bar G=k_FG\tens_s C(k,q)$
in general with cross relations $v\cdot x=s(x) x\cdot v$. On the
other hand, when $\sigma$ is an involutive automorphism this is
clearly a super tensor product of $\Z_2$-graded algebras.

\begin{lemma} When $\sigma$ is an involutive automorphism as in the
Clifford process, we have $\bar A=A\und\tens C(k,q)$ a super
tensor product. Moreover, applying twice with $q_1,q_2$ gives \[
\bar{\bar A}\isom A\und\tens C(k\oplus k,(q_1,q_2)). \]
\end{lemma}
\proof The super tensor product $A\und\tens C(k,q_1)$ contains
each factor as subalgebras with cross relations
$av=(-1)^{\xi(a)}va=\sigma(a)v$ where $\xi(a)$ is the $\Z_2$-degree
corresponding to $\sigma$. This is obviously the content of the
Clifford process. Applying twice we have
\[\bar{\bar A}=\bar
A\und\tens C(k,q_2)=(A\und\tens C(k,q_1))\und\tens
C(k,q_2)=A\und\tens(C(k,q_1)\und\tens C(k,q_2)\] using that the
super tensor product $\und\tens$ is an associative operation. We then
use Corollary~2.6 \eproof

This superalgebra periodicity can then be expressed in more usual
form using the following proposition.

\begin{propos} Let $\dim(V)=2m$ be even and $\sigma$
an involutive automorphism on $k_FG$ defined by $s$ of the form
$s(x)=(-1)^{\rho(x)}$ for a $\Z$-valued function $\rho$. Then
\[ k_F G\und\tens C(V,{\bf q})\isom k_{F'}G\tens C(V,{\bf q}),\]
where
\[ F'(x,y)= F(x,y) ((-1)^{m(2m-1)}q_1\cdots
q_{2m})^{\h(\rho(xy)-\rho(x)-\rho(y))}.\]
If $\imath=\sqrt{-1}\in k$ then $F'$ is cohomologous to $F$.
\end{propos}
\proof Here $\rho$ must differ from an additive character by an even
integer, hence $F'$ is well-defined. We
assume that $q_i=\pm 1$ so that $\mu\equiv (-1)^{m(2m-1)}q_1\cdots
q_{2m}=\pm 1$. If $\mu=1$ then $F'=F$. Otherwise if
$\imath=\sqrt{-1}\in k$
then clearly $F'=F\del s$ where $s(x)=\imath^{-\rho(x)}$. Here $\del
s(x,y)=s(x)s(y)/s(xy)$ is (an
exact) cocycle. The same calculation shows that
$(-1)^{\h(\rho(xy)-\rho(x)-\rho(y))}$ is a cocycle even if $\imath\notin
k$, so that $F'$ is necessarily a cocycle.
Now let $\phi:k_FG\tens C(V,{\bf q})\to k_{F'}G\tens
C(V,{\bf q})$ be defined by $\phi(x)=x (e_1\cdots
e_{2m})^{\rho(x)}$ when restricted to $k_FG$ and the identity on
$C(V,{\bf q})$. Here $\gamma\equiv
e_1\cdots e_{2m}$ implements the $\Z_2$-grading
automorphism of $C(V,{\bf q})$ as in Corollary~2.5 and has square
$\mu=\pm 1$ (these are in fact the only properties of $C(V,{\bf q})$
that we use, i.e. the same result applies for any superalgebra with
grading implemented by an element $\gamma$ with $\gamma^2=\pm 1$).
Then for all $x,y\in G$ we have
\[\phi(x\cdot_F y)=F(x,y)\phi(xy)=F(x,y)xy\gamma^{\rho(xy)}=
{F(x,y)\over F'(x,y)}x\cdot_{F'} y
\gamma^{\rho(x)+\rho(y)}\mu^{\h(\rho(xy)-\rho(x)-\rho(y))}
=x\gamma^{\rho(x)}y\gamma^{\rho(y)}\]
which is $\phi(x)\phi(y)$ as required. We also have
\[\phi(e_ix)=\phi(xe_i (-1)^{\rho(x)})=
x\gamma^{\rho(x)}e_i(-1)^{\rho(x)}=e_ix\gamma^{\rho(x)}
=\phi(e_i)\phi(x)\] as required. Hence $\phi$ (which is clearly
a linear isomorphism) is an algebra isomorphism as well. \eproof

The group and cocycle $F$ above are quite general but when we put in the
form in Section~2 for Clifford algebras we immediately obtain
\eqn{perlemm}{ C(V,{\bf q})\und\tens C(\pm)\isom C(V,-{\bf
q})\tens C(\pm).}
Here $C(k\oplus k,(q,q))$ is $C(2,0)=C(+)$ or $C(0,2)=C(-)$ and
$F'(x,y)=F(x,y)(-1)^{x\cdot y}$ in view of (\ref{rhoplus}), which
changes ${\bf q}$ to ${\bf -q}$. This along with the periodicity
Lemma~3.8 implies the usual periodicity properties for Clifford algebras
such as \eqn{per}{ C(0,n+2)\isom C(n,0)\tens \H,\quad
C(n+2,0)\isom C(0,n)\tens M_2(k).} The additional observation
$\H\tens\H\isom M_4(k)$ then gives the usual table with period 8
for all positive or all negative signatures (contructed
together). If $k$ is algebraically closed (or at least has
$\sqrt{-1}$) the situation is even simpler; we have periodicity 2
in the Clifford algebras i.e. \eqn{perC}{ C(2m)\isom
M_{2^m}(k),\quad C(2m+1)\isom M_{2^m}(k)\oplus M_{2^m}(k).}

Finally, it is possible to extend the Clifford process also to
representations.

\begin{propos} If $W$ is an irreducible representation of $A$
{\em not} isomorphic to
$W_\sigma$ defined by the action of $\sigma(a)$ then $\bar
W=W\oplus W_\sigma$ is an irreducible representation $\pi$ of
$\bar A$ obtained via the Clifford process with $q$. Here
\[\pi(v)=\left(\begin{matrix}{0&1\cr q&0}\end{matrix}\right), \quad
\pi(a)=\left(\begin{matrix}{a&0\cr
0&\sigma(a)}\end{matrix}\right)\] are the action on $W\oplus W$
in block form (here $\pi(a)$ is the explicit action of $a$ in the
direct sum representation $W\oplus W_\sigma$). If $W,W_\sigma$
{\em are} isomorphic then $W$ itself is an irreducible
representation of $\bar A$ for a suitable value of $q$.
\end{propos}
\proof Here $W_\sigma$ is the same vector space as $W$ but with
$a$ acting by $\sigma(a)$. Clearly we have a representation of
$\bar A=A\und\tens k[v]$ since
$\pi(v)\pi(a)=\pi(\sigma(a))\pi(v)$ and $\pi(v)^2=q$. If $U\subset
W\oplus W$ is a nonzero subrepresentation then it is also a
subrepresentation under $A$. The   projection to the first or
second part of the direct sum is $A$-equivariant hence its image
is either 0 or $W$ since $W,W_\sigma$ are irredicible. Hence a
non-zero $U$ has dimension at least that of $W$. If equal
dimension then one or other projection is an isomorphism of $U$
with $W$ or $W_\sigma$ but not both since these are not
isomorphic. But in this case the form of $\pi(v)$ implies a
contradiction. If greater dimension then consider the two maps
$W\to W\oplus W/U$ by embedding to each summand and quotienting.
The image has smaller dimension than $W$ hence both maps are zero
by $W$ irreducible. Hence $U=W\oplus W$.

If the two representaions $W$ and $W_\sigma$ are equivalent then
there exists an invertible linear map $\phi:W\to W$ such that
$\phi\rho(a)=\pi(\sigma(a))\phi$. We let $\pi(v)=\phi$. Note that
$\phi^2$ is central since $\sigma$ has order 2, hence $\phi^2=q$
for some $q\in k^*$. Thus $W$ itself extends to an irreducible
representation of $\bar A$. \eproof

Thus $k$ is an irreducible representation of $k=C(0,0)$
equivalent to its conjugate. Hence $k$ is also an irreducible
representation of $k[v]=C(1,0)$, the sign representation (say).
This is not equivalent to its conjugate under $\sigma$ (which is
the trivial representation). Hence $k^2$ is an irreducible
representation of $M_2(k)=C(2,0)$, its usual one, and so on. In
this way the natural representations over any field may be mapped
out.

\section{Spinor representations}

In this section we use the $k_FG$ method to obtain a new approach
to the spinor representations  for Clifford algebras. We have
already seen that Clifford algebras may be constructed as
braided-commutative algebras in a symmetric monoidal category,
where the braiding has the form (\ref{Psi}) in terms of a
$\Z_2^n$-grading. One may make several categorical constructions
along the lines of usual vector space constructions but with the
braiding. For example there is a braided tensor
product\cite{Ma:introm} algebra
$A\und\tens_\Psi A$ which acts on $A$ from the left and right.
The right action can be viewed as a left action using the
braided-commutativity of $A$. Here the braided tensor product and
the action are
\eqn{btens}{(a\tens b)(a'\tens b')=a\Psi(b\tens a')b'
=aa'\tens bb'(-1)^{\rho(b)\cdot \rho(a')+|b|\cdot|a'|},\quad
(a\tens b).c=abc,} when $A$ has braiding $\Psi$ of the form in
(\ref{Psi}) (so $|\ |$ is the $\Z_2^n$-grading and $\rho$ the
induced $\Z_2$-grading). The following is a variant of this
observation in which we work with the super tensor product
algebra and modify the action to compensate for this.

\begin{propos} Suppose that $\imath=\sqrt{-1}\in k$. If $A$ is a
$\Z_2^n$-graded braided-commutative algebra with respect to
$\Psi$ of the form in (\ref{Psi}) then the super tensor product
$A\und\tens A$ acts on $A$ by
\[ (a\tens b).c=abc (-1)^{|b|\cdot |c|}\imath^{\rho(b)}\]
where the $\rho(b)\equiv \sum_{i=1}^n |b|_i=|b|\cdot|b|$ is
viewed in $\Z$ rather than in $\Z_2$. Moreover, the action is a
$\Z_n$-graded and $\Z_2$-graded one.
\end{propos}
\proof Applying the action twice gives
\[ (a\tens b).((a'\tens b').c)=aba'b'c
(-1)^{|b'|\cdot|c|+|b|\cdot(|a'|+|b'|+|c|)}\imath^{\rho(b')+\rho(b)}\]
while the action of the super tensor product is
\[ (-1)^{\rho(b)\rho(a')}(aa'\tens
bb').c=aa'bb'c(-1)^{\rho(b)\rho(a')}(-1)^{(|b|+|b'|)\cdot|c|}
\imath^{\rho(b+b')}\]
which gives the same when we use braided-commutativity to write
$a'b=ba'(-1)^{\rho(a')\rho(b)+|a'|\cdot|b|}$ and when we note that
\[
\imath^{\rho(b+b')}=\imath^{\rho(b)}\imath^{\rho(b')}(-1)^{|b|\cdot|b'|}
\]
by writing $\rho(b)=|b|\cdot|b|$ (in other words the existence of
$\sqrt{-1}$ allows us to write the cocycle $(-1)^{|b|\cdot|b'|}$
as a group coboundary as in Proposition~3.9). Also, since the action
is given by the product in $A$ and this is $\Z_n$-graded it follows
that the representation here is a  $\Z_n$-graded and hence
$\Z_2$-graded one as well. \eproof

In particular, we can apply this result to any Clifford algebra
acting on itself. The super tensor product algebra is a Clifford
algebra on a vector space of twice the dimension by Corollary~2.6.

\begin{corol} If $\imath=\sqrt{-1}\in k$ then $C(V\oplus V,{\bf
q}\oplus{\bf q})\isom C(V,{\bf q})\und\otimes C(V,{\bf q})$ acts
on $C(V,{\bf q})$ by \[ (e_x\tens e_y).e_z=e_x\cdot e_y\cdot e_z
(-1)^{y\cdot
z}\imath^{\rho(y)}=e_{x+y+z}F(x,y)F(x+y,z)(-1)^{y\cdot
z}\imath^{\rho(y)}\] where $\rho(y)=y\cdot y\in \Z$. Moreover,
the action is irreducible and yields an isomorphism
\[ C(V\oplus V,{\bf q}\oplus {\bf q})\isom {\rm End}(C(V,{\bf q})).\]
\end{corol}
\proof  Here $C(V,{\bf q})$ is a superalgebra with the required
braided-commutativity from Corollary~2.3 as required. We write in
the explicit form of the product in the basis $\{e_x\}$ and the
additive group structure of $\Z_2^n$. This holds in fact for any
$k_F\Z_2^n$ algebra with $\CR$ of the required form. (Putting in
$F$ from Proposition~2.1 would give the action in the Clifford
algebra case even more explicitly.) For the irreducibility and
the identification with endomorphisms it suffices to show that
the action is faithfull (since the dimensions match). Thus
suppose that
\[0= \sum_{x,y\in\Z_2^n}c_{x,y}(e_x\tens e_y).e_z =c_{x,y}F(x,y)
F(x+y,z)(-1)^{y\cdot z}
\imath^{\rho(y)} e_{x+y+z}\] for all $z\in \Z_2^n$. Changing
variables to $x+y=x'$, as $x'$ varies the vectors $e_{x'+z}$ run
through a basis (since $\Z_2^n$ is a group). So
\[ 0=\sum_{y}c_{x'+y,y}F(x'+y,y)(-1)^{y\cdot z}\imath^{\rho(y)},
\quad\forall z,x'.\]
We dropped the $F(x',z)$ factor here since it is non-zero for all
$z,x'$. For each $x'$ fixed this is the $\Z_2^n$
Fourier-transform of a function of $y$, hence the funtion
vanishes for all $y$. Hence $c_{x,y}$ vanish, i.e. our action is
faithfull.

Also, the $\Z_2$-grading on the representation here is that of
$C(V,{\bf q})$ and coincides with the canonical one induced by
the action of the `top' element $e_{(1,\cdots ,1)}\tens
e_{(1,\cdots,1)}$ of $C(V\oplus V,{\bf q}\oplus{\bf q})$. From the
above it is given by
\[ (e_{(1,\cdots,1)}\tens e_{(1,\cdots 1)}).e_z= e_z F((1,\cdots,1),
(1,\cdots,1))(-1)^{\rho(z)}\imath^n=
\lambda (-1)^{\rho(z)} e_z,\quad \lambda=\imath^n
(-1)^{n(n-1)\over 2}\prod_i q_i.\] \eproof

Our action is in closed form, but the explicit action of the
generators on using the specific form of $F$ is \eqn{spinor}{
(e_i\tens
1).e_x=(-1)^{\sum_{j=1}^{i-1}x_j}q_i^{x_i}e_{x+(0,\cdots,1,\cdots
0)}, \quad (1\tens e_i).e_x=\imath
(-1)^{\sum_{j=1}^{i-1}x_j}(-q_i)^{x_i} e_{x+(0,\cdots,1,\cdots
0)},} where $(0,\cdots, 1,\cdots 0)$ denotes $1$ in the $i$'th
place. This construction is very different from but yields the
same result as the usual construction of spinors\cite{Gre:mul}
on the exterior algebra $\Lambda V$, which has the same dimension
as $C(V,{\bf q})$. Thus, we identify the bases $\{e_{i_1}\wedge
\cdots\wedge e_{i_p}\}$ and $\{e_{i_1}\cdots e_{i_p}\}$ in the two
cases. With a standard choice of polarisation (or complex structure) on
$V\oplus V$ the usual spinor representation (with $e_i\equiv
e_i\tens 1$ and $e_{n+i}\equiv 1\tens e_i$) is
\eqn{spinext}{ e_i.\omega=e_i\wedge\omega+i_{e_i}\omega,\quad
e_{n+i}.\omega=\imath(e_i\wedge \omega-i_{e_i}\omega),
\quad\forall \omega\in \Lambda V.} Here $i_v$ denotes the
interior product in $\Lambda V$ with the $\bf q$ norm. This
coincides with our action on any $e_x$ since, if $x_i=0$ only
$e_i\wedge$ contributes while if $x_i=1$ only $i_{e_i}$
contributes. The effect is therefore to change $x_i$ to $x_i+1$
mod 2 as in (\ref{spinor}). The coefficients also coincide.

Next, given any super representation of a super algebra $A$ one
has a usual representation of the cross product algebra
$A\lcross k\Z_2$ which is the bosonisation\cite{Ma:book} of $A$.
The latter is just defined by
adjoining a generator $v$ with $v^2=1$  and cross relations
$vav^{-1}=\sigma(a)$). An irreducible super representation $W$ of
$A$ extends to the bosonistion with $v.w=(-1)^{\rho(w)}w$ for all
$w\in W$, where $\rho$ is the super degree.

\begin{corol} When $\sqrt{-1}\in k$ the odd Clifford algebras
$C(V\oplus V\oplus k, ({\bf q}\oplus
{\bf q},q))$ are also represented irreducibly in the vector space
of $C(V,{\bf q})$. Here the action of $C(V\oplus V,{\bf
q}\oplus{\bf q})$ above is extended by the additional generator
$e_{2n+1}$ acting as $\lambda\sigma_W$ with $\lambda^2=q$.
\end{corol}
\proof Here $W=C(V,{\bf q})$ and $A=C(V\oplus V,{\bf q}\oplus{\bf
q})$.  The bosonisation consists in adjoining $v$ which clearly
equivalent to the Clifford process. By a minor rescaling of the
generator $v$ we adopt instead the relation $v^2=q$  for Clifford
process with parameter $q$ and identify it with $e_{2n+1}$ of
$C(V\oplus V\oplus k, ({\bf q}\oplus{\bf q},q))$. Note that $W$
remains irreducible since any submodule restricted to $C(V\oplus
V,{\bf q}\oplus{\bf q})$ must coincide with $W$. \eproof

One can also endow $A\lcross k \Z_2$ with a new super algebra
structure, with $v$ of degree 1. The extended representation $W$
is no longer a super representation
but one can be obtained by doubling it to $\bar W=W\oplus W$. Thus
$C(V,{\bf q})\oplus C(V,{\bf q})\isom C(V\oplus k,({\bf q},q))$
becomes a super representation of the odd Clifford algebra. In
this case applying the bosonisation again gives a representation
equivalent to the next higher even Clifford algebra
representation acting as in the above Corollary but acting on
$C(V\oplus k,({\bf q},q))$. One can also view these results from
the Clifford process point of view in the previous section.

As a very concrete example of our main result, consider the spinor
represenation for the Clifford algebra $C(0,4)$ in 4 Euclidean
dimensions. We work over $\C$
and by the above this can be considered as $\H\und\tens \H$ acting on
$\H$ where $\H=C(0,2)$ is the complex quaternions. Thus, a `Dirac
spinor' in physics is nothing other than an $\H$-valued
function. With basis $\{e_1,e_2\}$ of the 2-dimensional vector
space $V$ taken as the generators of $\H$, the spinor
action from Corollary~4.2 is
\eqn{spineucl}{( e_i\tens 1).\psi=e_i\psi,\quad
(1\tens e_i).\psi=\imath e_i \psi (-1)^{|\psi|_i}}
on a spinor of homogeneous degree $|\psi|\in \Z_2^2$. The right hand
side here uses the quaternion product. The construction of spinor
representations in terms of left and right actions of quaternions
have previously been alluded to in some contexts in the literature,
see for example\cite{Dix:div}. However, we are not aware of a general
treatment as above.

As an application, the standard Dirac operator on the
4-dimensional space $V\oplus V$ under the identification of
Corollary~2.6 is
\[ \dsl\psi=(e_1\tens 1)\del_1+(e_2\tens 1)\del_2+(1\tens
e_1)\del_3+(1\tens e_2)\del_4\]
where $\del_1$ is differentiation in the first basis direction of
$V\oplus V$, etc. So this becomes
\eqn{dirac}{ \dsl\psi=e_1(\del_1+\imath\del_3
(-1)^{|\psi|_1})\psi+e_2(\del_2+\imath\del_4(-1)^{|\psi|_2})\psi.}
It is possible to make this more explicit by taking a basis of $\H$,
namely
$1,e_1,e_2$ and $e_3\equiv e_1e_2$. Then
\[ e_ie_j=-\delta_{ij}+\eps_{ijk}e_k,\quad i,j,k=1,2,3\]
as usual, in terms of the Kronecker delta-function and the totally
antisymmetric tensor with $\eps_{123}=1$. We write a spinor as an
ordered pair $\psi=(\psi_0,\psi_i)$ with $i=1,2,3$
according to the components in
this basis. Finally, we write
\[ \nabla_1=\del_1+\imath\del_3,\quad \nabla_2=\del_2+\imath\del_4\]
and denote by $\bar \nabla$ the same expressions with $-\imath$.  Then
\[ (\dsl\psi)_0=-\bar\nabla_1\psi_1-\bar\nabla_2\psi_2,\quad
(\dsl\psi)_1=\nabla_1\psi_0+\bar\nabla_2\psi_3\]
\[ (\dsl\psi)_2=\nabla_2\psi_0-\bar\nabla_1\psi_3,\quad
(\dsl\psi)_3=\nabla_1\psi_2-\nabla_2\psi_1\]
using the relations in $\H$. If we define $\nabla_3=0$ and
$\vec\psi=(\psi_1,\psi_2,\bar\psi_3)$ then this can be written
compactly as
\eqn{diraccurl}{ (\dsl\psi)_0=-\bar\nabla\cdot\vec{\psi},\quad
\vec{\dsl\psi}=\nabla\psi_0+\bar\nabla\times \bar{\vec{\psi}}}
in terms of usual divergence, gradient and curl in 3 (complex)
dimensions and pointwise complex conjugation.


\begin{thebibliography}{10}


\bibitem{AlbMa:qua}
H. Albuquerque and S.~Majid.
\newblock Quasialgebra Structure of Octonions.
\newblock {\em J. Algebra}, 220:188--224, 1999.

\bibitem{Wen:cliff}
G.P.Wene.
\newblock A construction relating Clifford algebras and Cayley-Dickson
algebras.
\newblock {\em J.  Math. Phys.}, 25(8):2351--2353, 1984.

\bibitem{Lam:quad}
T.Y.Lam.
\newblock {\em The algebraic theory of quadratic forms}.
\newblock W.A.Benjamim, Inc. (Advanced Book Program), Massachusetts, 1973.


\bibitem{Dri:qua}
V.G. Drinfeld.
\newblock Quasi{H}opf algebras.
\newblock {\em Leningrad Math. J.}, 1:1419--1457, 1990.

\bibitem{Ma:book}
S.~Majid.
\newblock {\em Foundations of Quantum Group Theory}.
\newblock Cambridge Univeristy Press, 1995.

\bibitem{GurMa:bra}
D.I. Gurevich and S.~Majid.
\newblock Braided groups of {H}opf algebras obtained by twisting.
\newblock {\em Pac. J. Math.}, 162:27--44, 1994.


\bibitem{Mac:cat}
S.~Mac Lane.
\newblock {\em Categories for the Working Mathematician}.
\newblock Springer, 1974.
\newblock GTM vol. 5.



\bibitem{Ma:introm}
S.~Majid.
\newblock Algebras and {H}opf algebras in braided categories.
\newblock Volume 158 of {\em Lec. Notes in Pure and Appl. Math}, pages
55--105. Marcel Dekker, 1994.

\bibitem{Gre:mul}
W.H.~Greub.
\newblock {\em Multilinear Algebra.}
\newblock Graduate Texts in Maths, Vol. 136 (2nd edition).
Springer-Verlag, 1978.

\bibitem{Dix:div}
G.~Dixon.
\newblock {\em Division Algebras: Octonions, Quaternions,
Complex Numbers and the Algebraic Design of Physics}.
\newblock Kluwer, 1994.

\end{thebibliography}

\end{document}